\documentclass[12pt]{article}
\usepackage{amsmath,amsthm}
\usepackage{amssymb,latexsym}
\usepackage{enumerate}

\def\Q{{\mathbb Q}}
\def\Z{{\mathbb Z}}
\def\N{{\mathbb N}}

\def\O{{\cal O}}
\newtheorem{lemma}{Lemma}
\newtheorem{theorem}[lemma]{Theorem}

\title{
Relative power integral bases \\
in infinite families of quartic extensions \\
of quadratic fields 
}
\author{
Istv\'{a}n Ga\'{a}l 
\thanks{
        Research supported in part by K75566 and K100339 from the
        Hungarian National Foundation for Scientific Research
        }
\; and T\'\i mea Szab\'o 
\\
University of Debrecen, Mathematical Institute \\
            H--4010 Debrecen Pf.12., Hungary \\
            e--mail: igaal@science.unideb.hu, timi.taylor168@gmail.com  \\ \\
}

\begin{document}

\maketitle
\thispagestyle{empty}

\renewcommand{\thefootnote}{}



\renewcommand{\thefootnote}{\arabic{footnote}}
\setcounter{footnote}{0}

\begin{abstract}
\noindent 
We consider infinite parametric families of octic fields, that are quartic
extensions of quadratic fields. We describe all relative power integral
bases of the octic fields over the quadratic subfields.
\end{abstract}

\section{Introduction}
Let $K$ be an algebraic number field of degree $n$ with ring of integers $\O=\Z_K$.
The index of $\alpha\in\O$ (assumed $\alpha$ is a primitive element that is
$K=\Q(\alpha)$) is defined by
\[
I(\alpha)=(\O^+:\Z[\alpha]^+)
\]
that is the index of the additive group of $\Z[\alpha]$ in the additive group of $\O$.
The ring $\O$ is called {\it monogenic} if it is generated by a single element 
$\vartheta$ over $\Z$, that is $\O=\Z[\vartheta]$.
In this case the index of $\vartheta$, is equal to 1 and obviously
$1,\vartheta,\ldots,\vartheta^{n-1}$ is an integral basis of $\O$ called
{\it power integral basis}.
 
There is an extensive literature of monogenic number fields (cf. \cite{gaalbook})
involving also algorithms for determining all possible generators of power
integral bases. We also succeeded to determine all possible generators of
power integral bases also in some infinite parametric families of number fields,
see \cite{gpc5}, \cite{gaallettl}.

The index and the notion of power integral basis was extended also to the 
{\it relative case}, for relative extensions of number fields. Relative cubic extensions
were considered by I.Ga\'al and M.Pohst \cite{rel3} and relative quartic extensions by
I.Ga\'al and M.Pohst \cite{rel4} (cf. also \cite{gaalbook}).

Let $M$ be a number field and $K$ an extension field of $M$ of relative degree $n$.
Denote the rings of integers of $M$ and $K$ by $\Z_M$ and $\O=\Z_K$, respectively.
If $\alpha\in\O$ is a primitive element of $K$ over $M$ (that is $K=M(\alpha)$), then
the {\it relative index} of $\alpha$ over $M$ is defined by
\[
I_{K/M}(\alpha)=(\O^+:\Z_M[\alpha]^+),
\]
that is the index of the additive group of $\Z_M[\alpha]$ in the additive group of $\O$.
The relative index is equal to 1 if and only if $1,\alpha,\ldots,\alpha^{n-1}$
is a {\it relative power integral basis} of $\O$ over $\Z_M$.

The method of \cite{rel4} to determine relative power integral bases in
quartic relative extensions was a direct generalization of our method
\cite{gppsim} for quartic fields, however its application is much more complicated
technically.

\section{Results}

We are going to consider three infinite parametric families of octic fields $K=M(\xi)$
over their quadratic subfield $M$.
Our purpose is to describe the relative power integral bases 
of either $\O=\Z_K$ over $\Z_M$ (if the integer basis of $K$ is
known in a parametric form) or of $\O=\Z_M[\xi]$ over $\Z_M$
(otherwise). Note that in the later case $\xi$ itself is a generator
of a power integral basis but it is important to ask if there exist
any other generators of power integral bases.

The elements of $\O$ that only differ by unit factors 
and by translation by elements of $\Z_M$ are called {\it equivalent}.
Equivalent elements have the same relative indices.
There are only finitely many non-equivalent
generators of relative power integral bases.
We intend to determine all generators of power integral bases
up to equivalence.

\section{Results}

We shall consider three infinite families of relative quartic extensions
over quadratic fields.

\vspace{1cm}
\noindent
{\bf I.} Let $D>0$ be a square-free integer, $M=\Q(\sqrt{-D})$, 
$t\in\Z_M$ a parameter and let $\xi$ be a root of 
\begin{equation}
f(x)=x^4-t^2x^2+1\in \Z_M[x].
\label{f1}
\end{equation}
Let $K=M(\xi)$ and consider the relative power integral bases
of $\O=\Z_M[\xi]$ over $\Z_M$.

\begin{theorem}
For $|t|>245$ all non-equivalent generators of power integral bases of $\O$
over $\Z_M$ are given by
\[
\alpha=\xi,-t^2\xi+\xi^3,(1-t^4)\xi+t\xi^2+t^2\xi^3,(1-t^4)\xi-t\xi^2+t^2\xi^3,
t\xi^2+\xi^3,-t\xi^2+\xi^3.
\]
Moreover for $D=-3$ we also have
\[
\alpha=(1-\omega_3^2t)\xi+\omega_3\xi^2+\omega_3^2\xi^3
\]
with $\omega_3=(1+i\sqrt{3})/2$.
\label{th1}
\end{theorem}

\vspace{1cm}
\noindent
{\bf II.} Let $D>0$ be a square-free integer, $M=\Q(\sqrt{-D})$, 
$t\in\Z_M$ a parameter and let $\xi$ be a root of 
\begin{equation}
f(x)=x^4-4tx^3+(6t+2)x^2+4tx+1\in \Z_M[x].
\label{f2}
\end{equation}
Let $K=M(\xi)$ and consider the relative power integral bases
of $\O=\Z_M[\xi]$ over $\Z_M$.

\begin{theorem}
For $|t|>1544803$ all non-equivalent generators of power integral bases of $\O$
over $\Z_M$ are given by
\[
\alpha=\xi,(6t+2)\xi-4t\xi^2+\xi^3.
\]
\label{th2}
\end{theorem}

Note that the arguments to prove the above two theorems apply
results on the solutions of parametric relative Thue equations that
were only proved for large parameters. This makes necessary the assumptions on the
parameters.

\vspace{1cm}
\noindent
{\bf III.} Let $M=\Q(i)$, and let $t\in\Z_M$ be a parameter such that the polynomial
\begin{equation}
f(x)=x^4-itx^2+1\in \Z_M[x]
\label{f3}
\end{equation}
is irreducible over $M$. Denote by $\xi$ a root of $f(x)$. 

For this family of relative quartic extensions we shall separately
consider the Gaussian integer ($t\in\Z_M$) and rational integer ($t\in \Z$) 
parameters. The reason is that
if $t\in\N$ and $t^2+4$ is square free, then  the integral basis of $K$ is known
by B.K.Spearman and K.S.Williams \cite{sw}.

First consider Gaussian integer parameters. 
Let $t\in \Z_M\setminus \Z$ be a parameter such that $f(x)$ is irreducible
over $M$. Let $K=M(\xi)$,  set $\O=\Z_M[\xi]$
and consider the relative power integral bases of $\O$ over $\Z_M$.

\begin{theorem}
For the above parameters $t\in \Z_M\setminus \Z$ with $|t|<50$ 
all non-equivalent generators 
\[
\alpha=x\xi+y\xi^2+z\xi^3
\]
of power integral bases of $\O$ over $\Z_M$,
with $x,y,z\in\Z_M$, $\max(|x|,|y|,|z|)\leq 30$ are given by
\[
\alpha=\xi,-it\xi+\xi^3.
\]
\label{th3}
\end{theorem}

Next consider rational integer parameters.
By the results of B.K.Spearman and K.S.Williams \cite{sw}
for $t\in \N$, $t^2+4$ square free the polynomial $f(x)$ is irreducible
over $M$ and if $\xi$ is a root of $f(x)$, $K=M(\xi)$
then $\xi$ itself generates power integral basis over $M$.
Therefore we take $\O=\Z_K$ and
consider the relative power integral bases of $\O$ over $\Z_M$.

\begin{theorem}
For $t\in \N$, $t^2+4$ square free, $t\leq 100$,
all non-equivalent generators 
\[
\alpha=x\xi+y\xi^2+z\xi^3
\]
of power integral bases of $\O$ over $\Z_M$,
with $x,y,z\in\Z_M$, $\max(|x|,|y|,|z|)\leq 30$ are given by 
\[
\alpha=\xi,-it\xi+\xi^3.
\]
Moreover, for $t=1$ the following elements and their equivalents 
also generate power integral bases
of $\O$ over $\Z_M$:
\[
\alpha=3\xi+(1+i)\xi^2+2i\xi^3, 3\xi-(1+i)\xi^2+2i\xi^3, i\xi+(1+i)\xi^2+\xi^3,  i\xi-(1+i)\xi^2+\xi^3.
\]
\label{th4}
\end{theorem}

In the last two theorems we made direct calculation for the solutions of the
relative Thue equations involved. The solutions in general
are not known, that is why we only have results for small parameters.

\section{Preliminaries}

Our main tool throughout will be the application
the method of \cite{rel4} that we detail here.
This reduces the relative index form equation to a relative cubic equation
(Thue equation if irreducible) and some
relative quartic Thue equations. In our three families the cubic equation
will be reducible. To solve the quartic Thue equations we shall apply 
results on infinite parametric families of relative Thue equations.

Let $K$ be a quartic extension of the number field $M$ of degree $m$,
generated by a root
$\xi$ with relative minimal polynomial 
$f(x)=x^4+a_1x^3+a_2x^2+a_3x+a_4\in\Z_M[x]$.
Let $\O$ be either $\Z_K$ or $\Z_M[\xi]$.
We represent any $\alpha\in\O$ in the form
\begin{equation}
\alpha=\frac{1}{d}\left( a+x\xi+y\xi^2+z\xi^3  \right)
\label{alpha8}
\end{equation}
with coefficients 
$x,y,z\in\Z_M\; (1\leq i\leq 4)$ and with a common denominator $d\in\Z$.
Let $i_0=I_{K/M}(\xi)=(\O^+:\Z_M[\xi]^+)$,
\[
F(u,v)=u^3-a_2u^2v+(a_1a_3-4a_4)uv^2+(4a_2a_4-a_3^2-a_1^2a_4)v^3
\]
a binary cubic form over $\Z_M$ and
\begin{eqnarray*}
Q_1(x,y,z)&=&x^2-xya_1+y^2a_2+xz(a_1^2-2a_2)+yz(a_3-a_1a_2)+z^2(-a_1a_3+a_2^2+a_4)\\
Q_2(x,y,z)&=&y^2-xz-a_1yz+z^2a_2
\end{eqnarray*}
ternary quadratic forms over $\Z_M$.

\begin{lemma}(\cite{rel4})
If $\alpha$ of (\ref{alpha8}) satisfies
\[
I_{K/M}(\alpha)=1,
\]
then there is a solution $(u,v)\in \Z_M$ of
\begin{equation}
N_{M/\Q}(F(u,v))=\pm \frac{d^{6m}}{i_0}
\label{resolvent8}
\end{equation}
such that
\begin{eqnarray}
u&=&Q_1(x,y,z), \nonumber \\
v&=&Q_2(x,y,z).   \label{QV8}
\end{eqnarray}
\end{lemma}

For a given solution $u,v$ of (\ref{resolvent8}) we have to solve the system of
equations  (\ref{QV8}). For this purpose we use a method of 
L.J.Mordell \cite{mordell} to parametrize the solutions of the quadratic form equation
$Q_0(x,y,z)=uQ_2(x,y,z)-vQ_1(x,y,z)=0$ which is also explained in general
in \cite{rel4}. These detailes will be completely described in our proofs.

\section{Proofs}

\subsection{Proof of Theorem 1}

In our Lemma we substitute $a_1=0,a_2=-t^2,a_3=0,a_4=1$.
Equation (\ref{resolvent8}) is of the form
\[
F(u,v)=(u-2v)(u+2v)(u+t^2v)=\varepsilon
\]
where $\varepsilon$ is a unit in $M$. Therefore all factors of $F(u,v)$ must also
be units in $M$ which implies $v=0$ and $u$ a unit in $M$.
The equation $Q_0=0$ implies $Q_2=0$, that is
\begin{equation}
y^2-xz-z^2t^2=0.
\label{q01}
\end{equation}
A non-trivial solution of it is $x_0=-t^2,y_0=0,z_0=1$. Hence 
the solutions of (\ref{q01}) can be parametrized in the form 
\begin{eqnarray}
x&=&-t^2r+p\nonumber\\
y&=&q\label{xyz11}\\
z&=&r\nonumber
\end{eqnarray}
with $p,q,r\in M,r\neq 0$ (cf. \cite{mordell}, \cite{rel4}). Substituting these
expressions into (\ref{q01}) we obtain $rp=q^2$. We multiply all equations of 
(\ref{xyz11}) by $p$ and replace $rp$ by $q^2$, whence
\begin{eqnarray}
kx&=&-t^2q^2+p^2\nonumber\\
ky&=&pq\label{xyz12}\\
kz&=&q^2\nonumber
\end{eqnarray}
with some $k\in M$.
In case $\Z_M$ has unique factorization, we
multiply these equations by the square of the common denominators of
$p,q$ and divide them by the square of the gcd of $p,q$. This way 
we can replace the parameters $k,p,q$ by parameters in $\Z_M$. 
(Remark that if $\Z_M$ has no unique factorization, then principally the same 
argument is followed involving ideals of $\Z_M$. The detailed procedure can
be found in \cite{rel4}.)

Considering the coefficients of $p^2,pq,q^2$ in the
representation of $kx,ky,kz$ in (\ref{xyz12}) and calculating the determinant of the
corresponding 3x3 matrix we confer that $k$ must also be a unit (cf. \cite{rel4}).
Finally, substituting the representation (\ref{xyz12}) into $Q_1(x,y,z)=u$ we obtain
\[
p^4-t^2p^2q^2+q^4=k^2u.
\]
(The second equation $Q_2(x,y,z)=v$ vanishes on both sides.)

Applying the results of Theorem 2 of V.Ziegler \cite{ziegler} we can describe
the solutions of this relative quartic Thue equation, hence also $x,y,z$ in
(\ref{xyz12}) which gives the result of our Theorem 1. Note that 
Theorem 2 of V.Ziegler \cite{ziegler} is only valid for $|t|>245$.

\subsection{Proof of Theorem 2}

In our Lemma we substitute $a_1=-4t,a_2=6t+2,a_3=4t,a_4=1$.
Equation (\ref{resolvent8}) is of the form
\[
F(u,v)=(u+2v)(u-(2-2t)v)(u-(2+8t)v)=\varepsilon
\]
where $\varepsilon$ is a unit in $M$. Therefore all factors of $F(u,v)$ must also
be units in $M$ which implies $v=0$ and $u$ a unit in $M$.
Again $Q_0=0$ implies $Q_2=0$, that is
\begin{equation}
y^2-xz+4tyz+(6t+2)z^2=0.
\label{q02}
\end{equation}
A non-trivial solution of it is $x_0=6t+2,y_0=0,z_0=1$. Hence 
the solutions of (\ref{q02}) can be parametrized in the form 
\begin{eqnarray}
x&=&(6t+2)r+p\nonumber\\
y&=&q\label{xyz21}\\
z&=&r\nonumber
\end{eqnarray}
with $p,q,r\in M,r\neq 0$. Substituting these
expressions into (\ref{q02}) we obtain $q^2=r(p-4tq)$. We multiply all equations of 
(\ref{xyz21}) by $p-4tq$ and replace $r(p-4tq)$ by $q^2$, whence
\begin{eqnarray}
kx&=&p^2-4tpq+(6t+2)q^2\nonumber\\
ky&=&pq-4tq^2\label{xyz22}\\
kz&=&q^2\nonumber
\end{eqnarray}
with some $k\in M$.
Similarly as in I we replace the parameters
$k,p,q$ by parameters in $\Z_M$. 
Again we confer that $k$ must also be a unit.
Finally, substituting the representation (\ref{xyz22}) into $Q_1(x,y,z)=u$ we obtain
\[
p^4-4tp^3q+(6t+2)p^2q^2+4tpq^3+q^4=k^2u.
\]
(The second equation $Q_2(x,y,z)=v$ again vanishes on both sides.)

Applying the results of Theorem 2 of B.Jadrijevi\v c and V.Ziegler
\cite{jadrziegler} we can describe
the solutions of this relative quartic Thue equation, hence also $x,y,z$ in
(\ref{xyz22}) which gives the result of our Theorem 2. Note that 
Theorem 2 of B.Jadrijevi\v c and V.Ziegler \cite{jadrziegler}
is only valid for $|t|>1544803$.

\subsection{Proofs of Theorem 3 and Theorem 4}

In this proof we use the result of B.K.Spearman and K.S.Williams \cite{sw}.
They consider octic fields generated by a root of the polynomial
$x^8+(t^2+2)x^4+1=(x^4-itx^2+1)(x^4+itx^2+1)$. 
We set $M=\Q(i)$, $t\in\Z_M$ a parameter such that
$f(x)=x^4-itx^2+1$ is irreducible over $M$. Let $\xi$ be a root of $f(x)$
and $K=M(\xi)$. 

We take $t\in\Z_M\setminus \Z$, $\O=\Z_M[\xi]$ for the purposes of Theorem 3.
However, if $t\in \N$, $t^2+4$ squarefree, the result of
B.K.Spearman and K.S.Williams \cite{sw} implies, that $\xi$ itself
generates power integral basis of $K$ over $M$, that is why in Theorem 4
we take $\O=\Z_K$. 

In the following we consider the two cases (according to
Theorem 3 and Theorem 4, respectively) together.
In our Lemma we substitute $a_1=0,a_2=-it,a_3=0,a_4=1$.
Equation (\ref{resolvent8}) is of the form
\[
F(u,v)=(u+itv)(u-2v)(u+2v)=\varepsilon
\]
where $\varepsilon$ is a unit in $M$. Therefore all factors of $F(u,v)$ must also
be units in $M$ which implies $v=0$ and $u$ a unit in $M$.
By $Q_0=0$ we again get $Q_2=0$, that is
\begin{equation}
y^2-xz-itz^2=0.
\label{q03}
\end{equation}
A non-trivial solution of it is $x_0=1,y_0=0,z_0=0$. Hence 
the solutions of (\ref{q03}) can be parametrized in the form 
\begin{eqnarray}
x&=&r\nonumber\\
y&=&p\label{xyz31}\\
z&=&q\nonumber
\end{eqnarray}
with $p,q,r\in M,r\neq 0$. Substituting these
expressions into (\ref{q03}) we obtain $rq=p^2-itq^2$. We multiply all equations of 
(\ref{xyz31}) by $q$ and replace $rq$ by $p^2-itq^2$, whence
\begin{eqnarray}
kx&=&p^2-itq^2\nonumber\\
ky&=&pq\label{xyz32}\\
kz&=&q^2\nonumber
\end{eqnarray}
with some $k\in M$.
Similarly as before we replace the parameters
$k,p,q$ by parameters in $\Z_M$. Again we confer that $k$ must also be a unit.
Finally, substituting the representation (\ref{xyz32}) into $Q_1(x,y,z)=u$ we obtain
\begin{equation}
p^4-itp^2q^2+q^4=k^2u.
\label{t4}
\end{equation}
(The second equation $Q_2(x,y,z)=v$ again vanishes on both sides.)

Unfortunately the solutions of this parametric relative 
Thue equation are not known in general. We could only test "small"
solutions of the equation for "small" parameters $t$.

For the case of Theorem 3 we let $t\in\Z_M\setminus \Z$ run through
the range $t\leq 50$ (such that $f(x)$ is irreducible over $M$),
and tested the solutions of equation (\ref{t4})
for $|q|^2\leq 1000$ with the possible right hand sides $\pm 1,\pm i$.
(The range $|q|^2\leq 1000$ covers the range $\max(|x|,|y|,|z|)\leq 30$
by (\ref{xyz32}).)
For given $t,q$ and right hand side we calculated the roots $p$ of the
equation and tested if it is in $\Z_M$. In this range we only found trivial solutions,
that is $(p,q)=(1,0),(0,1)$ and their associates. Up to unit factors these solutions imply 
$(x,y,z)=(1,0,0),(-it,0,1)$ by (\ref{xyz32}).

For the case of Theorem 4 we let $t\in N$ run thought the values from $1\leq t\leq 100$ with
$t^2+4$ squarefree. For $|q|^2\leq 1000$ and the possible right hand sides $\pm 1,\pm i$
we calculated the corresponding roots of equation (\ref{t4}) and tested if it is in
$\Z_M$. 
(The range $|q|^2\leq 1000$ again corresponds to the range $\max(|x|,|y|,|z|)\leq 30$
by (\ref{xyz32}).)
In addition to the trivial solutions $(p,q)=(1,0),(0,1)$ we only found additional
solutions for $t=1$, namely $(p,q)=(1,1+i),(-1,1+i),(1+i,1),(1+i,-1)$ which give
the solutions listed in Theorem 4 by (\ref{xyz32}).

\vspace{1cm}
\noindent
{\it Remark 1.} There are usually several posible parametrization of the solutions
of the quadratic form equation $Q_0(x,y,z)=0$. It is not always straitforward to 
find a suitable parametrization leading to a known family of
parametric Thue equations. 

\noindent
{\it Remark 2.} The direct calculations involved in Theorems 3 and 4 were
performed in Maple on a simple PC. The CPU time took a couple of hours.

\noindent
{\it Remark 3.} Our experience supports the conjecture that the assertions of
Theorems 3 and 4 are true in general (not only for small parameters $t$ and
small values of $x,y,z$).

\end{document}